\documentclass{amsart}
\usepackage{latexsym,amsmath,amssymb}
\usepackage{cite}
\newtheorem{thm}{Theorem}[section]
\newtheorem{cor}[thm]{Corollary}
\newtheorem{exam}[thm]{Example}
\newtheorem{lem}[thm]{Lemma}
\newtheorem{prop}[thm]{Proposition}
\theoremstyle{definition}\newtheorem{defn}[thm]{Definition}
\theoremstyle{remark}
\newtheorem{rem}[thm]{Remark}
\numberwithin{equation}{section}

\begin{document}

\title[]
{Continuity of conditional expectation in Orlicz spaces}

\author{\sc\bf A. Hosseini and Y. Estaremi}
\address{\sc Y. Estaremi}
\email{y.estaremi@gu.ac.ir}
\address{Department of Mathematics, Faculty of Sciences, Golestan University, Gorgan, Iran.}
\address{\sc A. Hosseini}
\email{hoseiniseyedali299@gmail.com}
\address{Department of Mathematics, Faculty of Sciences, Golestan University, Gorgan, Iran.}
\email{}
\address{}
\thanks{}

\thanks{}



\date{}

\dedicatory{}

\commby{}

\begin{abstract}
	
	The continuity of conditional expectation on Orlicz spaces is investigated. Indeed, we provide some necessary and sufficient conditions on a sequence $\{\mathcal{A}_n\}_{n\in\mathbb{N}}$ of $\sigma$-subalgebras for $L^{\varphi}$-convergence of the related conditional expectations. Our results generalize
	similar results in $L^p$-spaces.
\end{abstract}

\maketitle

\section{ \sc\bf Introduction}
The study of conditional expectations and their convergence properties has evolved alongside the development of modern probability theory. Also the convergence of conditional expectations on increasing $\sigma$-algebras is a cornerstone of stochastic analysis, with roots in Kolmogorov’s axiomatization and Doob’s martingale theory.
The study of $L^p$-convergence of conditional expectations for increasing $\sigma$-algebras was significantly advanced by Donald Burkholder in the 1960s–1970s and Claude Dellacherie and Paul-André Meyer in the 1970s-1980s. Their work unified martingale theory with functional analysis and stochastic processes. For more information about the continuity of conditional expectation one can see \cite{7,8,4} and references therein. In \cite{7}, the authors provide a necessary and sufficient condition on a sequence $\{\mathcal{A}_n\}_{n\in\mathbb{N}}$ of $\sigma$-subalgebras that
assures $L^p$-convergence of the conditional expectations. Their results generalized
the $L^p$-martingales, the Fetter and the Boylan (equi-convergence) theorems. In this paper we are going to extend the results of \cite{7} to the setting of Orlicz spaces.\\

The continuous convex function  $\varphi:\mathbb{R}\rightarrow\mathbb{R}$ is called a Young function,  whenever 
\begin{itemize}
	\item $\varphi(x)=0$ if and only if $x=0$.
	\item $\varphi$ is even, i.e.,
	$\varphi(-x)=\varphi(x)$, for all
	$x\in\mathbb{R}$.
	\item $\lim_{x\rightarrow\infty}\frac{\varphi(x)}{x}=\infty,  \lim_{x\rightarrow\infty}\varphi(x)=\infty$. 
\end{itemize}
With each Young function $\varphi$, one can associate another convex function $\psi:\mathbb{R}\rightarrow\mathbb{R}^{+}$ having similar
properties, which is defined by   
$$
\psi(y):=\sup\{x|y|-\varphi(x):x\geq0\}, \quad y\in\mathbb{R}
$$
The $\psi$ is called the complementary Young fonction of $\varphi$ and $(\varphi, \psi)$ is called a pair of complementary Young function.

We say $\varphi$ satisfies the $\Delta_2$-condition, if there exist constants $k>0$ and $x_0>0$, such that $\varphi(2x)\leq k\varphi(x)$, for every $x\geq x_0$. We denotes this condition by $\varphi\in\Delta_2$.

Let $(\Omega, \mathcal{A}, \mu)$ be a finite measure space and $\mathfrak{D}\subseteq \mathcal{A}$ be a 
 $\sigma-$subalgebra. For a finite  $\sigma-$subalgebra $\mathfrak{D}\subseteq \mathcal{A}$, the conditional expectation operator associated with $\mathfrak{D}$ is the mapping 
$f\rightarrow E\left(f\mid\mathfrak{D}\right)$, 
defined for all non-negative,  measurable functions $f$ as well as for all $f\in L^1(\mathcal{A})$ and $f\in L^{\infty}(\mathcal{A})$, where $E\left(f\mid\mathfrak{D}\right)$, by the Radon-Nikodym theorem, is the unique $\mathcal{A}$-measurable function satisfying 
$$
\int_{D}f d\mu= \int_{D}E\left(f\mid\mathfrak{D}\right) d\mu, \quad \forall D\in\mathfrak{D}. 
$$ 
The function $E\left(f\mid\mathfrak{D}\right)$  is called conditional expectation of $f$ with respect to $\mathfrak{D}$. 
As an operator on $L^1(\mathcal{A})$ and $L^{\infty}(\mathcal{A})$, $E\left(\cdot\mid\mathcal{A}\right)$ is idempotent and $E\left(L^{\infty}(\mathcal{A}) \mid \mathfrak{D}\right)= L^{\infty}(\mathfrak{D})$ and $E\left(L^{1}(\mathcal{A}) \mid \mathfrak{D}\right)= L^{1}(\mathfrak{D})$. Thus it can be defined on all interpolation spaces of $L^1$ and $L^{\infty}$ such as Orlicz spaces \cite{9}. 
Let $\varphi$ be a Young function and $L^0(\mathcal{A})$ be the space of all $\mathcal{A}$-measurable real values functions on $\Omega$.    Then the set of $\mathcal{A}$-measurable functions  
$$
L^{\varphi}(\mathcal{A})
=L^{\varphi}(\Omega, \mathcal{A}, \mu):= \Big\{f\in L^0(\mathcal{A}):\exists  \alpha>0, \int_{\Omega}\varphi(\alpha f) d\mu<\infty \Big\}, 
$$
is a Banach space, with respect to the norm  $N_{\varphi}\left(.\right)$ that is defined as:
$$
N_{\varphi}\left(f\right)=\|f\|_{\varphi}:=\inf\{k>0:\ \int_{\Omega}\varphi\Bigl(\frac{|f|}{k}\Bigr) d\mu\leq1\}, 
$$ 
for all $f\in L^{\varphi}(\mathcal{A})$.  
The pair $\Bigl(L^{\varphi}(\Omega, \mathcal{A}, \mu), N_{\varphi}(.)\Bigr)$  is called an Orlicz space \cite{6}. 
For more details about Orlicz spaces see \cite{5,6}. 

Since a conditonal expectation depends on two variables, natural questions to ask are the following: 
\begin{enumerate}
\item [a] If $\{f_n\}_{n=1}^{\infty}\subseteq L^1(\mathcal{A})$ is a sequence which converges to $f\in L^1(\mathcal{A})$, when will  $\{E\left(f_n\mid\mathfrak{D}\right)\}_{n\in\mathbb{N}}$ converge to $E\left(f\mid\mathfrak{D}\right)$?
\item [b] If $\{\mathcal{A}_n\}_{n\in\mathbb{N}}$ is a sequence of $\sigma$-subalgebra of $\mathcal{A}$ which  converges to a  $\sigma$-subalgebra $\mathfrak{D}$ of $\mathcal{A}$ in some sense, when will  $\{E\left(f\mid\mathcal{A}_n\right)\}_{n\in\mathbb{N}}$ converge? 
\item [c] If the sequence   $\{f_n\}_{n\in\mathbb{N}}$  converges to $f$ and $\{\mathcal{A}_n\}_{n\in\mathbb{N}}$ converge to $\mathfrak{D}$, when will $\{E\left(f_n\mid\mathcal{A}_n\right)\}_{n\in\mathbb{N}}$ converge? 
\end{enumerate} 
For case (a), several answers are known, for instance, the continuity of conditional expectations notably established by using dominated convergence theorem in \cite{1}.  
Fore cases (b) and (c), if the sequence $\{\mathcal{A}_n\}_{n\in\mathbb{N}}$ is increasing (decreasing), i.e., 
if $\mathcal{A}_n\subseteq\mathcal{A}_{n+1}$ and $\mathfrak{D}=\bigvee_{n=1}^\infty\mathcal{A}_n$ ($\mathcal{A}_{n+1}\subseteq\mathcal{A}_n$ and $\mathfrak{D}=\bigcap_{n=1}^\infty\mathcal{A}_n$) for all $n\in\mathbb{N}$, then $\{E\left(f\mid\mathcal{A}_n\right)\}_{n\in\mathbb{N}}$ converges to $E\left(f\mid\mathfrak{D}\right)$ and if $\{f_n\_{n\in\mathbb{N}}\}$ converges to $f$ in $L^p(\mathcal{A})$, then $\{E\left(f_n\mid\mathcal{A}_n\right)\}_{n\in\mathbb{N}}$ converges to $E\left(f\mid\mathfrak{D}\right)$ in $L^p$, for $1\leq p<\infty$. These results are known as the martingale convergence theorem \cite{1,2,3}. Fetter in \cite{4} have shown that if $\{\mathcal{A}_n\}_{n\in\mathbb{N}}$ is a sequence of $\sigma$-subalgebra of $\mathcal{A}$ such that $\mathcal{A}=\bigvee_{n=1}^\infty\bigcap_{m=n}^{\infty}\mathcal{A}_m=\bigcap+{n=1}^{\infty}\bigvee_{m=n}^\infty\mathcal{A}_m$, then $E\left(f\mid\mathcal{A}_n\right)$ converges to $E\left(f\mid\mathcal{A}\right)$ in $L^p$, for $f\in L^p(\mathcal{A})$ and $p\geq1$.
 But $E(\left(f\mid\mathcal{A}_n\right))$ does not necessarily converge in $L^{\infty}$, for $f\in L^{\infty}\left(\mathcal{A}\right)$.  
Some necessary and sufficient conditions for the continuity of $E\left(f\mid\mathcal{A}_n\right)$ with respect to $L^p$-convergence are provided in \cite{7}. 

The $ L^p $ spaces are well-established for their simplicity and wide applications in analysis and probability theory [1]. However, Orlicz spaces, by introducing a Young function $\varphi$, provide greater flexibility in modeling nonlinear and asymmetric behaviors, which are crucial for addressing complex problems in harmonic analysis, statistics, and beyond. This paper generalizes the results from [1] in $L^p$-spaces to Orlicz spaces, tackling challenges posed by the nonlinear properties of $\varphi$, thus extending the applicability of these results to a broader class of function spaces. 

A classical question, that is originally studied in $L^p$-spaces, for $1\leq p<\infty$, is to determine under what conditions the sequence $\{E\left(f\mid \mathcal{A}_n\right)\}_{n\in\mathbb{N}}$ converges in $L^p(\mathcal{A})$. 

In this paper, we assume that 
$
(\Omega, \mathcal{A}, \mu)$ is a complete finite measure space and provide some necessary and sufficient conditions for the continuity of the sequence $\{E\left(f\mid\mathcal{A}_n\right)\}_{n\in\mathbb{N}}$ in the Orlicz space $L^{\varphi}(\Omega, \mathcal{A}, \mu)$. 

\section{ \sc \bf $\mu$-Convergence}  
In this section the $\mu$-Convergence of a sequence of $\sigma$-subalgebras $\{\mathcal{A}_n\}_{n\in\mathbb{N}}$ will be defined and the relation between the $\mu$-Convergence of $\{\mathcal{A}_n\}_{n\in\mathbb{N}}$ and the convergence of the sequence of corresponding conditional expectations will be characterized. Now first we recall the definition of $\mu$-Convergence of a sequence of $\sigma$-subalgebras of a finite $\sigma$-algebra.
\begin{defn}\cite{7}
	Let $(\Omega, \mathcal{A}, \mu)$ be a measure space,   $\{\mathcal{A}_n\}_{n\in\mathbb{N}}$ be a sequence of $\sigma$-subalgebras of $\mathcal{A}$ and  $\mathfrak{D}$ be a $\sigma$-subalgebra of $\mathcal{A}$.  We say that $\{\mathcal{A}_n\}_{n \in \mathbb{N}}$  $\mu$-converges to $\mathfrak{D}$, and denote it by $\mathcal{A}_n \xrightarrow{\mu} \mathfrak{D}$, if for any $D\in\mathfrak{D}$ there exists a sequence $\{A_n: A_n\in\mathcal{A}_n\}_{n \in \mathbb{N}}$  such that $ \lim_{n \rightarrow \infty} \mu(A_n \Delta D)= 0$. 
\end{defn}
In $\mu$-Convergence the limit is not unique. For example, $\{\emptyset, \Omega\}$ is the smallest $\sigma$-subalgebra of $\mathcal{A}$ such that every  $\{\mathcal{A}_n\}_{n\in\mathbb{N}}$ $\mu$-converges to it.  We will focus to finding the largest  $\sigma$-subalgebra in that $\{\mathcal{A}_n\}_{n\in\mathbb{N}}$ $\mu$-approaches to it. In the sequel for every sequence $\{\mathcal{A}_n\}_{n \in \mathbb{N}}$ of $\sigma$-subalgebras of $\mathcal{A}$ we define 
$$ \mathcal{A}_\mu := \{A \in \mathcal{A} : \exists A_n \in \mathcal{A}_n \ \ \ s.t.  \ \ \  \lim_{n \rightarrow \infty} \mu(A_n \Delta A) = 0\}. $$
The collection of measurable sets $\mathcal{A}_{\mu}$ is a complete $\sigma$-subalgebra of $\mathcal{A}$. Moreover, we use the symbole  $\bigvee_{n=1}^{\infty}\mathcal{A}_n$ for the smallest $\sigma$-subalgebra of $\mathcal{A}$ containing $\bigcup_{n=1}^{\infty}\mathcal{A}_n$.  We set  $\overline{\mathcal{A}}:=\bigcap^{\infty }_{m=1}{\bigvee^{\infty }_{n=m}{{\mathcal{A}}_n}}$ and  $\underline{\mathcal{A}}=\bigvee^{\infty }_{m=1}{\bigcap^{\infty }_{n=m}{{\mathcal{A}}_n}}$,  and call them  upper and lower limits of the sequence $\{\mathcal{A}_n\}_{n=1}^{\infty}$, respectively.  Now, we recall some lemmas from \cite{7} for later use.
\begin{lem}\label{lem1}\cite{7}
	Let $(\Omega, \mathcal{A}, \mu)$ be a measure space and  $\{\mathcal{A}_n\}_{n \in \mathbb{N}}$ be a sequence of $\sigma$-subalgebras of $\mathcal{A}$ and $\mathfrak{D}$  be a $\sigma$-subalgebra of $\mathcal{A}$. Then $\mathcal{A}_n \xrightarrow{\mu} \mathfrak{D}$ if and only if $\mathfrak{D} \subseteq A_\mu$.
\end{lem}
For every sequence of  $\sigma$-subalgebras  $\{\mathcal{A}_n\}_{n\in\mathbb{N}}$ of $\mathcal{A}$, the collection  $\mathcal{A}_{\mu}$ is the largest $\sigma$-subalgebra of $\mathcal{A}$ such that  $\mathcal{A}_n\xrightarrow{\mu}\mathcal{A}_{\mu}$. 
\begin{lem}\cite{7}
	Suppose that $\left(\Omega, \mathcal{A}, \mu\right)$ is a measure space and  $\{\mathcal{A}_n\}_{n \in \mathbb{N}}$ is a sequence of $\sigma$-subalgebras of $\mathcal{A}$. Then 
	$$ \underline{A} \subseteq A_\mu \subseteq \overline{A}. $$
\end{lem}
In the next lemma, we find an equivalent condition to the $\mu$-convergence of a sequence of $\sigma$-subalgebras in the Orlicz space.
\begin{lem}\label{lem3}
	Let $\left(\Omega, \mathcal{A}, \mu\right)$ be a measure space, $\{\mathcal{A}_n\}_{n \in \mathbb{N}}$ be a sequence of $\sigma$-subalgebras of $\mathcal{A}$ and also $\mathfrak{D}$  be a $\sigma$-subalgebra of $\mathcal{A}$. Moreover, let  $\varphi$ be a Young function. Then  $\mathcal{A}_n \xrightarrow{\mu} \mathfrak{D}$ if and only if for every $D \in \mathfrak{D}$, there exists a sequence $\{A_n: A_n\in\mathcal{A}_n\}_{n\in\mathbb{N}}$ such that the corresponding sequence  $\{|\chi_{A_n} - \chi_D|\}_{n \in \mathbb{N}}$ converges to zero in the Orlicz space $\Bigl(L^\varphi(\mathcal{A}),   N_{\varphi}\left(\cdot\right)\Bigr)$. 
	\begin{proof}  
		Let $\mathcal{A}_n \xrightarrow{\mu} \mathfrak{D}$. Then for every $D \in \mathfrak{D}$ there exists a sequence $\{A_n: A_n\in\mathcal{A}_n\}_{n\in \mathbb{N}}$ such that $\lim_{n \rightarrow \infty} \mu(A_n \Delta D) = 0$. Now fix  $m \in \mathbb{N}$ and set $\gamma_m := \frac{1}{m} > 0$, and   $f_n := \chi_{A_n} - \chi_D$, for every $n \in \mathbb{N}$. So we have $|f_n(x)| = 1$, for  $x \in A_n \Delta D$ and $|f_n(x)|= 0$, for  $x \in\Omega\setminus\left(A_n \Delta D\right)$. Hence
		\begin{align*}
			\int_\Omega \varphi\left(\frac{|f_n(x)|}{\gamma_m}\right) d\mu &= \int_{A_n \Delta D} \varphi\left(\frac{1}{\gamma_m}\right) d\mu\\
			&= \mu(A_n \Delta D) \ \varphi(m).
		\end{align*}
		Since $\lim_{n \rightarrow \infty} \mu(A_n \Delta D)= 0$, then for every $\varepsilon > 0$ there exists $N_{\varepsilon, m} \in \mathbb{N}$ such that 
		$$\mu(A_n \Delta D) < \frac{\varepsilon}{\left(\varepsilon + 1\right) \varphi(m)}, $$
		for all $n>N_{\varepsilon, m}$ and consequently 
		$$ \int_\Omega \varphi\left(\frac{|f_n(x)|}{\gamma_m}\right) d\mu < \frac{\varepsilon}{\left(\varepsilon + 1\right) \varphi(m)} \varphi(m) = \frac{\varepsilon}{\varepsilon + 1} < 1, $$
		for every $n > N_{\varepsilon, m}$. This implies that for every $\varepsilon > 0$ and every $m \in \mathbb{N}$, there exists $N_0 := N_{\varepsilon, m}$ such that
		\begin{align*}
			N_\varphi(f_n) \leq \gamma_m,\tag{1}
		\end{align*}
		for all $n \geq N_0$. Here,  we suppose on the contrary that the sequence  $\{f_n\}_{n=1}^\infty$ does not converge to zero in $\Bigl(L^\varphi(\mathcal{A}), N_\varphi\Bigr)$, i.e.,   
		$$ \exists \varepsilon_0 > 0\quad \forall N \in \mathbb{N}\quad \exists n_0 \geq N:\quad N_\varphi(f_{n_0}) \geq \varepsilon_0. $$
		Let $m_0 \in \mathbb{N}$ such that $ \gamma_{m_0}= \frac{1}{m_0} < \varepsilon_0$. Then
		$$ \exists \varepsilon_0 > 0 \ \  \exists m_0 \in \mathbb{N} \ \ \forall N \in \mathbb{N} \ \ \exists n_0 \geq N: \ \ N_\varphi(f_{n_0}) \geq \varepsilon_0 > \frac{1}{m_0} = \gamma_{m_0}.$$
		This contradicts the inequality (1). Therefore, we get that  the sequence $\{f_n\}_{n=1}^\infty$ and so $\{|f_n|\}_{n=1}^\infty$ is  converges to zero in $\Bigl(L^\varphi(\mathcal{A}), N_\varphi\Bigr)$.\\
		Suppose that for every $D \in \mathfrak{D}$, there exists a sequence $\{A_n: A_n\in\mathcal{A}_n\}_{n \in \mathbb{N}}$ such that the sequence $\{|\chi_{A_n} - \chi_D|\}_{n=1}^\infty$ converges to zero in $\Bigl(L^\varphi(\mathcal{A}), N_\varphi\Bigr)$. This implies that the sequence $\{N_\varphi(\chi_{A_n} - \chi_D)\}_{n=1}^\infty$ converges to zero in $\mathbb{R}$ and so by the definition of norm $N_{\varphi}(.)$,  there exists a sequence of positive numbers $\{\lambda_n\}_{n=1}^\infty$,  converging to $0$,  such that
		\begin{align*}
			\varphi\left(\frac{1}{\lambda_n}\right) \mu(A_n \Delta D) \leq 1,\tag{2} 
		\end{align*}
		for all $n\in\mathbb{N}$. If we set $f_n := \chi_{A_n} - \chi_D$, then $f_n\in L^{\varphi}$, for all $n \in \mathbb{N}$. Now, if  $\lim_{n \rightarrow \infty} \mu(A_n \Delta D) \neq 0$,  Then
		$$ \exists \varepsilon_0 > 0\quad \forall n \in \mathbb{N}\quad \exists k_n > n:\quad \mu(A_{k_n} \Delta D) \geq \varepsilon_0. $$
		Since $\lambda_n \rightarrow 0$, then we have $\frac{1}{\lambda_n} \rightarrow \infty$ and so $\varphi\left(\frac{1}{\lambda_{k_n}}\right) \rightarrow \infty$ as $n \rightarrow \infty$. Consequently,
		$$ \int_\Omega \varphi\left(\frac{|f_{k_n}|}{\lambda_{k_n}}\right) d\mu = \varphi\left(\frac{1}{\lambda_{k_n}}\right) \mu(A_{k_n} \Delta D)> 1, $$ 
		for sufficiently large $n$, which contradicts (2).	Hence we get the result. 
	\end{proof}
\end{lem}  
It is known that for each $f \in L^\varphi(\mathcal{A})$ we have 
$$ \int_\Omega \varphi\left(\frac{f}{N_\varphi(f)}\right) d\mu \leq 1.$$
In the next lemma we recall the relation between convergence in the Orlicz space and integral convergence. 
\begin{lem}\cite{6}\label{lem5}
	Let $(\Omega, \mathcal{A}, \mu)$ be a measure space, $\varphi$ be a Young function and	$f \in L^\varphi(\Omega,\mathcal{A},\mu)$. If $\lim_{n \rightarrow \infty} N_\varphi(f_n-f) = 0$, for a sequence $\{f_n\}_{n=1}^\infty\subseteq L^\varphi(\mathcal{A})$, then $$\lim_{n\rightarrow\infty} \int_{\Omega} \varphi(f_n - f)  d\mu=0.$$
\end{lem}
It is well know that for finite measure space $(\Omega, \mathcal{A}, \mu)$, convex function $\varphi:\mathbb{R}\rightarrow\mathbb{R}$ and measurable function $f:\Omega \rightarrow \mathbb{R}$ such that $\int_{\Omega} \varphi(f) d\mu$ and $\int_{\Omega} f d\mu$ exist, we have 
$$\varphi\left(\frac{1}{\mu(\Omega)} \int_{\Omega} f  d\mu\right) \leq \frac{1}{\mu(\Omega)} \int_{\Omega} \varphi(f)  d\mu.$$ 
It is known as Generalized Jensen's Inequality. In the following theorem, we will find a relationship between the $\mu$-convergence of a sequence of $\sigma$-subalgebras and the convergence of their corresponding  conditional expectations in the Orlicz space $L^\varphi (\mathcal{A})$. 
\begin{thm}
	Let $(\Omega, \mathcal{A}, \mu)$ be a measure space, $\varphi$ be a Young function,  $\{\mathcal{A}_n\}_{n \in \mathbb{N}}$ be a sequence of $\sigma$-subalgebras of $\mathcal{A}$ and $\mathfrak{D}$ be a $\sigma$-subalgebra of $\mathcal{A}$. Then the following hold.\\ 
	\textbf{(a)} If $\mathcal{A}_n \xrightarrow{\mu} \mathfrak{D}$ and $\varphi \in \Delta_2$, then $E\left(f\mid \mathcal{A}_n\right) \xrightarrow{L^\varphi(\mathcal{A})} E\left(f \mid \mathfrak{D}\right) = f$,  for all $f \in L^\varphi(\mathfrak{D})$.
	\\ 
	\textbf{(b)} If $E\left(f\mid \mathcal{A}_n\right) \xrightarrow{L^\varphi(\mathcal{A})} E\left(f\mid \mathfrak{D}\right)= f$, for all $f \in L^\varphi(\mathfrak{D})$, then $\mathcal{A}_n \xrightarrow{\mu} \mathfrak{D}$.
	\begin{proof}
		\textbf{(a)} 
		Let $\mathcal{A}_n \xrightarrow{\mu} \mathfrak{D}$ and $D \in \mathfrak{D}$. Then there  exists a sequence $\{A_n: A_n\in\mathcal{A}_n\}_{n\in\mathbb{N}}$ such that $\lim_{n \rightarrow \infty} \mu(A_n \Delta D) = 0$. For $k>0$  and $k_1=k_2=\frac{k}{2}$, we have 
		\begin{align*} 
			\int_{\Omega}\varphi \Bigl(\frac{
				E\left({\chi_{D}}\mid{\mathcal{A}_n}\right)-\chi_{D}}{k}\Bigr)& d\mu \\
			=&\int_{\Omega}\varphi \Bigl(\frac{E\left({\chi_{D}}\mid{\mathcal{A}}_n\right)-\chi_{A_n}  +{\chi_{A_n}}-\chi_{D}}{k}\Bigr) d\mu\\ 
			=&\int_{\Omega}\varphi \Bigl(\frac{E\left({\chi_{D}}\mid{\mathcal{A}}_n\right)-{\chi_{A_n}}}{k} + \frac{{\chi}_{A_n}-{\chi}_D}{k}\Bigr) d\mu\\
			=&\int_{\Omega}\varphi \Bigl(\frac{E\left({\chi_{D}}\mid{\mathcal{A}_n}\right) -E\left({\chi_{A_n}}\mid{\mathcal{A}_n}\right)}{k} +  \frac{\chi_{A_n}-{\chi_{D}}}{k}\Bigr) d\mu \\
			=& \int_{\Omega}\varphi \Bigl(\frac{E\left({\chi_{D}-\chi_{A_n}}\mid{\mathcal{A}_n}\right)}{k}+ \frac{\chi_{A_n}-\chi_{D}}{k}\Bigl) d\mu\\
			=& \int_{\Omega} \varphi\Bigl(\frac{E\left({\chi_{D}-\chi_{A_n}}\mid{\mathcal{A}_n}\right)}{k_1}\cdot \frac{k_1}{k}+\frac{\chi_{A_n}-\chi_{D}}{k_2}\cdot  \frac{k_2}{k}\Bigr) d\mu\\
			\le&  \frac{k_1}{k}  \int_{\Omega} \varphi  \Big(\frac{E\left({\chi_{D}-\chi_{A_n}}\mid{\mathcal{A}_n}\right)}{k_1}\Bigr) d\mu\\
			&  +\frac{k_2}{k}\int_{\Omega}\varphi \Bigl(\frac{\chi_{A_n}-\chi_{D}}{k_2}\Bigr) d\mu\ \\
			\leq& \int_{\Omega}\varphi \Bigl( \frac{E\left({\chi_{D}-\chi_{A_n}}\mid{\mathcal{A}_n}\right)}{k_1}\Bigr) d \mu\\
			& +\int_{\Omega}\varphi \Big(\frac{\chi_{A_n}-\chi_{D}}{k_2}\Bigl) d\mu \\
			\leq&  \int_{\Omega}E\Bigl(\varphi(\frac{\chi_{D}-\chi_{A_n}}{k_1})\mid\mathcal{A}_n\Bigr) d\mu\\
			& +\int_{\Omega}\varphi\Bigl(\frac{\chi_{A_n}-\chi_{D}}{k_2}\Bigr) d\mu \\
			=&\int_{\Omega}\varphi\Bigl(\frac{\chi_{D}-\chi_{A_n}}{k_1}\Bigr) d\mu +\int_{\Omega}\varphi\Bigl(\frac{\chi_{A_n}-\chi_{D}}{k_2}\Bigr) d\mu\\
			\leq& \int_{\Omega}\varphi\Bigl(\frac{|\chi_{D}-\chi_{A_n}|}{k_1}\Bigr) d\mu +\int_{\Omega}\varphi\Bigl(\frac{|\chi_{A_n}-\chi_{D}|}{k_2}\Bigr) d\mu\\
			=&\mu(A_n\Delta D) \ \varphi(\frac{1}{k_1})+\mu(A_n\Delta D) \ \varphi(\frac{1}{k_2})\\
			=&2 \mu(A_n\Delta D) \ \varphi(\frac{2}{k}).
		\end{align*}
		Hence for any $k>0$, we have
		$$ \int_{\Omega}\varphi \Bigl(\frac{E\left({\chi_{D}}\mid{\mathcal{A}_n}\right)-\chi_{D}}{k}\Bigr) d\mu\leq 2 \mu(A_n\Delta D) \ \varphi(\frac{2}{k}),$$
		and so  
		$$ \Big\{k>0:  2\mu(A_n\Delta D)\cdot \varphi(\frac{2}{k})\leq1\Big\}\subseteq \Big\{k>0:  \int_{\Omega}\varphi \Bigl(\frac{E\left({\chi_{D}}\mid{\mathcal{A}_n}\right)-\chi_{D}}{k}\Big) d\mu\leq1\Big\}. $$
		Taking infimum over $k>0$ on both sides, we get that 
		\begin{flalign*}
			N_{\varphi}(E\left({\chi_{D}}\mid{\mathcal{A}_n}\right)-\chi_{D})&\leq \inf\Big\{k>0: 2\mu(A_n\Delta D)   \cdot \varphi(\frac{2}{k})\leq 1 \Big\}\\
			&=\inf\Big\{k>0: k\geq\frac{2}{\varphi^{-1}(\frac{1}{2\mu(A_n\Delta D)})}\Big\}\\
			&=\frac{2}{\varphi^{-1}(\frac{1}{2\mu(A_n\Delta D)})}.
		\end{flalign*}
		Since $\lim_{n \rightarrow \infty} \mu(A_n \Delta D) = 0$, then we have
		\begin{align*}
			N_{\varphi}(E\left({\chi_{D}}\mid{\mathcal{A}_n}\right)-\chi_{D})\xrightarrow{n\rightarrow\infty}0.\tag{1}
		\end{align*}
		So for the characteristic functions we have the result. For the next step we consider simple function i.e., we let $\psi = \sum_{i=1}^M a_i \chi_{D_i}\in L^{\varphi}(\mathfrak{D})$ be a simple  function, in which, $D_i\in\mathfrak{D}$, $a_i\in\mathbb{R}$ and $M\in\mathbb{N}$. Since $E$ is a bounded linear operator, then by the relation (1) we have 
		\[
		E\left(\psi\mid \mathcal{A}_n\right) \xrightarrow{L^\varphi(\mathcal{A})} E\left(\psi\mid \mathfrak{D}\right)= \psi.  \tag{2}
		\] 
		Now, we consider $f \in L^\varphi(\mathfrak{D})$. Since $\varphi \in \Delta_2$, then the set of $\mathfrak{D}$-measurable simple functions is dense in $L^\varphi(\mathfrak{D})$. So  there exists a sequence $\{\psi_m\}_{m=1}^\infty \subseteq L^\varphi(\mathfrak{D})$ of simple functions such that
		\[
		\psi_m \xrightarrow{L^\varphi(\mathfrak{D})} f. \tag{3}
		\]
		Again by linearity of conditional expectation,  
		for any $N \in \mathbb{N}$, we have
		\begin{align*}
			N_{\varphi}(E\left({f}\mid{\mathcal{A}_n}\right)-f)&\le N_{\varphi}(E\left({f}\mid{\mathcal{A}_n}\right)-E\left({\psi_N}\mid{\mathcal{A}_n}\right))\\
			& \ \ \ +N_{\varphi}(E\left({\psi_N}\mid{\mathcal{A}_n}\right)-\psi_N)\\
			& \ \ \ +N_{\varphi}({\psi}_N -f)\\
			&=N_{\varphi}(E\left({f-\psi_N}\mid{\mathcal{A}_n}\right))\\
			&\ \ \ +N_{\varphi}(E\left({\psi_N}\mid{\mathcal{A}_n}\right)-\psi_N)\\
			& \ \ \ +N_{\varphi}({\psi}_N -f)\\
			&\le 2N_{\varphi}(\psi_N-f)+N_{\varphi}(E\left({\psi_N}\mid{\mathcal{A}_n}\right)-\psi_N).
		\end{align*}
		Here first we take limit on $n$ and then take limit on $N$, when $n, N$$\rightarrow$$\infty$, and also by using (2) and (3), we get 
		\[
		\lim_{n \rightarrow \infty} N_{\varphi}(E\left(f \mid \mathcal{A}_n\right) - f) = 0.
		\]
		Therefore  $E\left(f \mid \mathcal{A}_n\right) \xrightarrow{L^\varphi(\mathcal{A})} f$, which proves the first part.\\
		\textbf{(b)} Let  $E\left(f \mid \mathcal{A}_n\right) \xrightarrow{L^\varphi(\mathcal{A})} E\left(f\mid \mathfrak{D}\right)= f$, for all $f \in L^\varphi(\Omega, \mathfrak{D}, \mu)$. Fix $D \in \mathfrak{D}$. Since $\chi_D \in L^\varphi(\Omega, \mathfrak{D}, \mu)$, then $E\left(\chi_D \mid \mathcal{A}_n\right) \xrightarrow{L^\varphi(\mathcal{A})} E\left(\chi_D\mid \mathfrak{D}\right) = \chi_D$. If we set $A_n := \{x \in \Omega : E\left(\chi_D \mid \mathcal{A}_n\right)(x) > \frac{1}{2}\}$, for every $n \in  \mathbb{N}$, then we have 
		$$\mu(A_n \Delta D) = \mu(A_n \cap D^c) + \mu(A_n^c \cap D)$$
		and  $A_n, A_n^c\in\mathcal{A}_n$, for all $n \in \mathbb{N}$. Also  
		\begin{align*}
			\frac{1}{2}\mu (A_n\Delta D) &=\frac{1}{2}\int_{\Omega}{\chi_{A_n}}{\chi_{D^c}} d\mu +\frac{1}{2}\int_{\Omega}\chi_{{A_n^c}}\chi_{D} d\mu\\
			&\le \int_{\Omega}\chi_{A_n}\chi_{D^c} |E\left({\chi_{D}}\mid{\mathcal{A}_n}\right)| d\mu +\int_{\Omega}\chi_{A_n^c}\chi_{D} |1-E\left({\chi_{D}}\mid{\mathcal{A}_n}\right)| d\mu\\ 
			&=\int_{\Omega}\chi_{A_n}\chi_{D^c} |E\left({\chi_{D}}\mid{\mathcal{A}_n}\right)-\chi_{D}| d\mu +\int_{\Omega}\chi_{A_n^c}\chi_{D} |E\left({\chi_{D}}\mid{\mathcal{A}_n}\right)-\chi_{D}| d\mu\\
			&\le \int_{\Omega}\chi_{A_n} |E\left({\chi_{D}}\mid{\mathcal{A}_n}\right)-\chi_{D}| d\mu +\int_{\Omega}\chi_{A_n^c} |E\left({\chi_{D}}\mid{\mathcal{A}_n}\right)-\chi_{D}| d\mu\\ 
			&= \int_{\Omega} |E\left({\chi_{D}}\mid{\mathcal{A}_n}\right)-\chi_{D}| d\mu. 
		\end{align*} 
		This implies that
		$$ \lim_{n \rightarrow \infty} \frac{1}{2} \mu(A_n \Delta D) \leq \lim_{n \rightarrow \infty} \int_\Omega |E\left(\chi_D \mid \mathcal{A}_n\right) - \chi_D| d\mu. $$
		Now, we show that $\lim_{n \rightarrow \infty} \int_\Omega |E\left(\chi_D \mid \mathcal{A}_n\right) - \chi_D| d\mu = 0$.  To this end, set $f_n := |E\left(\chi_D \mid \mathcal{A}_n\right) - \chi_D|$, for all $n \in \mathbb{N}$. By assumption $N_\varphi(f_n) \rightarrow 0$ and  the Lemma \eqref{lem5}, we have $\lim_{n \rightarrow \infty} \int_\Omega \varphi(f_n) d\mu = 0$. Hence $\lim_{n \rightarrow \infty} \frac{1}{\mu(\Omega)}\int_\Omega \varphi(f_n) d\mu = 0$ and so 
		$$ \varphi^{-1}\Bigl(\frac{1}{\mu(\Omega)}\int_{\Omega}\varphi(f_n)  d\mu\Bigr)\rightarrow0.$$
		Now by Generalized Jensen's Inequality, we get that $\lim_{n}\int_{\Omega}f_n  d\mu=0.$ This  completes the proof.
	\end{proof}
\end{thm}
Here in the next lemma we investigate the relation between convergence of a sequence of $\sigma$-subalgebras to their upper limit and the convergence of the sequence of corresponding conditional expectations.
\begin{lem}
	Let $(\Omega, \mathcal{A}, \mu)$ be a measure space, $\varphi\in \Delta_2$ be a Young function  and  $\{\mathcal{A}_n\}_{n=1}^\infty$ be a sequence of $\sigma$-subalgebras of $\mathcal{A}$. Then $\mathcal{A}_n \xrightarrow{\mu} \overline{\mathcal{A}}$ if and only if 
	$$ E\left(f\mid \mathcal{A}_n\right) \xrightarrow{L^\varphi(\mathcal{A})} E\left(f \mid \overline{\mathcal{A}}\right), $$
	for all $f \in L^\varphi(\mathcal{A})$. 
	\begin{proof}
		For $f \in L^\varphi(\mathcal{A})$, there exists a scaler $\alpha$ such that $\int_\Omega \varphi(\alpha |f|) d\mu < \infty$. Clearly, we get that 
		$$\int_\Omega \varphi(\alpha |E\left(f \mid \mathcal{A}_n\right)|) d\mu \leq \int_\Omega E\left({\varphi(\alpha |f |)}\mid{\mathcal{A}_n}\right) d\mu = \int_\Omega \varphi(\alpha |f|) d\mu< \infty.$$
		$(\Rightarrow):$ Assume that  $\mathcal{A}_n \xrightarrow{\mu} \overline{\mathcal{A}}$. By definition we have the inclusion   $\mathcal{A}_{\mu} \subseteq \overline{\mathcal{A}}$ and also by  Lemma \eqref{lem1}  we have $\overline{\mathcal{A}}\subseteq \mathcal{A}_{\mu}$. So $\overline{\mathcal{A}}= \mathcal{A}_{\mu}$. Since  $\underline{\mathcal{A}}\subseteq \mathcal{A}_{\mu}\subseteq  \mathcal{A}_{\perp}\subseteq\overline{\mathcal{A}}$, then $\mathcal{A}_{\mu}=  \mathcal{A}_{\perp}=\overline{\mathcal{A}}$. Hence by Theorem \eqref{thm4.5}  we get  that 
		$$E\left({f}\mid{\mathcal{A}_n}\right)\xrightarrow{L^{\varphi}(\mathcal{A})}E\left({f}\mid{\overline{\mathcal{A}}}\right),$$
		for all $f \in L^\varphi(\mathcal{A})$.\\
		($\Leftarrow$): 
		Suppose that  
		$$E\left(f \mid \mathcal{A}_n\right) \xrightarrow{L^\varphi(\mathcal{A})} E\left(f \mid \overline{\mathcal{A}}\right),$$
		for all $f \in L^\varphi(\mathcal{A})$. We need to show that $\mathcal{A}_n \xrightarrow{\mu} \overline{\mathcal{A}}$, which follows directly from Theorem \eqref{thm4.5}. 
	\end{proof}
\end{lem}

\section{ \sc \bf $\perp$-Convergence} 
In this section the $\perp$-Convergence of a sequence of $\sigma$-subalgebras $\{\mathcal{A}_n\}_{n\in\mathbb{N}}$ will be defined and the relation between the $\perp$-Convergence of $\{\mathcal{A}_n\}_{n\in\mathbb{N}}$ and the convergence of the sequence of corresponding conditional expectations will be characterized. \\

Let $(\Omega, \mathcal{A}, \mu)$ be a measure space, $\varphi$ be a Young function,  $\mathfrak{D}$ be a  $\sigma$-subalgebra of $\mathcal{A}$  and $f \in L^\varphi(\mathcal{A})$. The orthogonal complement of the conditional expectation of $f$ with respect to $\mathfrak{D}$, is denoted by $E_{\mathfrak{D}}^\perp(f)$ and is  defined as  
$$E_{\mathfrak{D}}^\perp(f) := f - E\left(f\mid \mathfrak{D}\right).$$
Now first we recall the definition of $\perp$-Convergence of a sequence of $\sigma$-subalgebras of a finite $\sigma$-algebra.
\begin{defn}
	Let $(\Omega, \mathcal{A}, \mu)$ be a measure space,  $\varphi\in \Delta_2$ be a Young function,  $\{\mathcal{A}_n\}_{n=1}^\infty$ be a sequence of $\sigma$-subalgebras of $\mathcal{A}$ and $\mathfrak{D}$ be another  $\sigma$-subalgebra of $\mathcal{A}$.  We say that  $\{\mathcal{A}_n\}_{n=1}^\infty$ is \textbf{Orthogonally convergent} to $\mathfrak{D}$ in $L^\varphi(\Omega, \mathcal{A}, \mu)$,  denote $\mathcal{A}_n \xrightarrow{\perp} \mathfrak{D}$, if for all sequences $\{A_n : A_n \in \mathcal{A}_n\}_{n=1}^\infty$ we have
	$$E_{\mathfrak{D}}^\perp(\chi_{A_n}) \xrightarrow{weakly} 0, \quad \text{in } L^\varphi(\Omega, \mathcal{A}, \mu).$$
\end{defn}

\begin{rem} 
	If $(\varphi, \psi)$ is a pair of complementary Young's functions and $\varphi\in \Delta_2$,  then $(L^\varphi)^* = L^\psi$ \cite{6}. Hence $E_{\mathfrak{D}}^\perp(\chi_{A_n}) \xrightarrow{weakly} 0$ if and only if
	$$\int_\Omega E_{\mathfrak{D}}^\perp(\chi_{A_n}) g d\mu \rightarrow 0,$$
	for every $g \in L^\psi(\Omega, \mathcal{A}, \mu)$.
\end{rem}

\begin{rem}
	For a sequence $\{\mathcal{A}_n\}_{n \in \mathbb{N}}$ of $\sigma$-subalgebras of $\mathcal{A}$, we can find more than one $\sigma$-subalgebra like $\mathfrak{D}$ of $\mathcal{A}$ such that $\mathcal{A}_n\xrightarrow{\perp}\mathfrak{D}$. This means that in the orthogonally convergence the limit is not unique. For example,  $\mathcal{A}$ is the largest $\sigma$-subalgebra such that every  $\{\mathcal{A}_n\}_{n\in\mathbb{N}}$, $\perp$-approaches to $\mathcal{A}$. We will focus to finding the smallest $\sigma$-subalgebra that $\{\mathcal{A}_n\}_{n\in\mathbb{N}}$ $\perp$-approaches to it. Now,  we define 
	$$W_\varphi := \left\{g \in L^\varphi(\mathcal{A}) : \exists A_{n_k} \in \mathcal{A}_{n_k} \ \ \ s.t. \ \ \  \chi_{A_{n_k}} \xrightarrow{weakly} g \right\}. $$
	Let $\mathcal{A}_{\perp}$ be the smallest complete  $\sigma$-subalgebra of $\mathcal{A}$ that the elements of $W_{\varphi}$ are  $\mathcal{A}_{\perp}$-measureable. \end{rem}
In the following Lemma we find an equivalent conditions to orthogonally convergence of the sequence $\{\mathcal{A}_n\}_{n=1}^\infty$ to $\mathfrak{D}$.
\begin{lem}\label{lem10}
	Let $(\Omega, \mathcal{A}, \mu)$ be a measure space, $\varphi\in\Delta_2$ be  a Young function,    $\{\mathcal{A}_n\}_{n=1}^\infty$  be a sequence of $\sigma$-subalgebras of $\mathcal{A}$ and $\mathfrak{D}$ be an arbitrary $\sigma$-subalgebra of $\mathcal{A}$. Then $\mathcal{A}_n \xrightarrow{\perp} \mathfrak{D}$ if and only if $\mathcal{A}_\perp \subseteq \mathfrak{D}$.
	\begin{proof}
		Suppose that $\mathcal{A}_n \xrightarrow{\perp} \mathfrak{D}$ and $g\in W_{\varphi}$.  This implies that there exists a  sequence $\{A_{n_k} : A_{n_k} \in \mathcal{A}_{n_k}\}_{k=1}^\infty$ such that  $\chi_{A_{n_k}}\xrightarrow{weakly}g$ and so  for every $h \in L^\psi(\mathcal{A})$ we have 
		\begin{align*}
			\int_{\Omega}h E\left({g}\mid{\mathfrak{D}}\right)d\mu &=\int_{\Omega}  E\left({h}\mid{\mathfrak{D}}\right) g  d\mu \\
			&=\lim_{k} \int_{\Omega}  E\left({h}\mid{\mathfrak{D}}\right) \chi_{A_{n_k}} d\mu \\ 
			&=\lim_{k} \int_{\Omega} h  E\left({\chi_{A_{n_k}}}\mid{\mathfrak{D}}\right) d\mu\\ 
			&=\lim_k \Bigl(\int_{\Omega} h (E\left({\chi_{A_{n_k}}}\mid{\mathfrak{D}}\right)-\chi_{A_{n_k}})   d\mu  +\int_{\Omega} h  \chi_{A_{n_k}} d\mu \Bigr)\\
			&=\lim_{k}\int_{\Omega} h\ {\chi}_{A_{n_k}}  d\mu\\
			&=\int_{\Omega}h g d\mu.
		\end{align*} 
		Hence $g$ and $E\left({g}\mid{\mathfrak{D}}\right)$ are almost every where equal and so $g$ is $\mathfrak{D}$-measurable function. Therefore $W_\varphi$ consists of $\mathfrak{D}$-measurable functions.  Consequently by definition of $\mathcal{A}_\perp$, we have $\mathcal{A}_\perp\subseteq\mathfrak{D}$.\\ 
		Conversely, assume that $\mathcal{A}_\perp \subseteq \mathfrak{D}$, but $\mathcal{A}_n \not\xrightarrow{\perp} \mathfrak{D}$. So there exists a  sequence $\{A_n : A_n \in \mathcal{A}_n\}_{n \in \mathbb{N}}$ such that $E_{\mathfrak{D}}^\perp(\chi_{A_n}) \not\xrightarrow{weakly} 0$. Hence, there exists  $g \in L^\psi(\mathcal{A})$  such that $\int_\Omega E_{\mathfrak{D}}^\perp(\chi_{A_n})\ g d\mu \not\rightarrow 0$ and so there exists  $\varepsilon > 0$ and a subsequence $\{A_{n_k}\}_{k \in \mathbb{N}}$ such that $|\int_\Omega E_{\mathfrak{D}}^\perp(\chi_{A_{n_k}}) g d\mu| \geq \varepsilon$, for all $k \in \mathbb{N}$. Since $\{\chi_{A_{n_k}}\}_{k \in \mathbb{N}}$  is a bounded sequence in $L^{\infty}(\mathcal{A})$ ($\|\chi_{A_{n_k}}\|_{\infty}=1$, for all $k\in\mathbb{N}$), then by the Banach–Alaoglu Theorem and the fact that $L^{\infty}(\mathcal{A})\cong L^1(\mathcal{A})^*$, it has a subsequence that converges in the weak* topology to some function $h\in L^{\infty}(\mathcal{A})\subseteq L^{\varphi}(\mathcal{A})$. Without loss of generality, we can assume $\chi_{A_{n_k}}\xrightarrow{weak^*}h$. Since the underlying measure space is finite, then $L^{\psi}(\mathcal{A})\subseteq L^{1}(\mathcal{A})$, and so the weak* convergence implies that 
		$$\lim_{k\rightarrow \infty}\int_{\Omega}\chi_{A_{n_k}}g d\mu =\int_{\Omega}hgd\mu,$$
		for all $g\in L^{\psi}(\mathcal{A})$. This is exactly the weak convergence  $\chi_{A_{n_k}} \xrightarrow{weakly} h$, in $L^{\varphi}(\mathcal{A})$.  So $h \in W_\varphi$. Therefore $h$ is $\mathcal{A}_\perp$-measurable. Since $\mathcal{A}_\perp \subseteq \mathfrak{D}$, then $h$ is $\mathfrak{D}$-measurable. Therefore
		\begin{align*}
			\lim_{k}\int_{\Omega}g E^{\perp}_{\mathfrak{D}}({{\chi }_{{\mathcal{A}}_{n_k}}}) d\mu &=\lim_{k}\int_{\Omega} \Bigl( g \chi_{{\mathcal{A}}_{n_k}}-g  E\left({{\chi }_{{\mathcal{A}}_{n_k}}}\mid{\mathfrak{D}}\right)\Bigr)d\mu\\
			&=\lim_{k}\int_{\Omega} \Bigl( g {\chi}_{{\mathcal{A}}_{n_k}}-E\left({g}\mid{\mathfrak{D}}\right) {\chi }_{{\mathcal{A}}_{n_k}}\Bigr)d\mu\\
			&=\int_{\Omega}\Bigl(g h-E\left({g}\mid{\mathfrak{D}}\right) h\Bigr) d\mu\\
			&=\int_{\Omega} \Bigl(g h-g E\left({h}\mid{\mathfrak{D}}\right)\Bigr) d\mu\\
			&=\int_{\Omega} \bigl(g h-g h\bigr)d\mu \\
			&=0.
		\end{align*}

		This is a contradiction, because   $|\int_\Omega E_{\mathfrak{D}}^\perp(\chi_{A_{n_k}}) g d\mu| \geq \varepsilon$, for each 
		$k\in\mathbb{N}$. 
	\end{proof}
\end{lem}
In the sequel we prove the inclusions $\underline{\mathcal{A}}\subseteq\mathcal{A}_\perp \subseteq \overline{\mathcal{A}}$ in the setting of Orlicz spaces. 
\begin{lem}
	Let $(\Omega, \mathcal{A}, \mu)$ be  a measure space, $\varphi\in \Delta_2$ be a Young function and $\{\mathcal{A}_n\}_{n \in \mathbb{N}}$ be a sequence of $\sigma$-subalgebras of $\mathcal{A}$. Then $\underline{\mathcal{A}}\subseteq\mathcal{A}_\perp \subseteq \overline{\mathcal{A}}$.
	\begin{proof}
		First we prove the inclusion  $\underline{\mathcal{A}}\subseteq\mathcal{A}_\perp$. To this end, fix $m \in \mathbb{N}$ and assume $A \in \bigcap_{n=m}^\infty \mathcal{A}_n$. Hence we can put $A_n := A$, for every $n \geq m$ and so $\chi_A \in W_\varphi$. Thus $\chi_A$ is $\mathcal{A}_\perp$-measurable. As we have chosen $A \in \bigcap_{n=m}^\infty \mathcal{A}_n$ arbitrary, we have $\bigcap_{n=m}^\infty \mathcal{A}_n \subseteq \mathcal{A}_\perp$ and since $\mathcal{A}_\perp$ is a $\sigma$-algebra, we get that $\bigcup_{m=1}^\infty \bigcap_{n=m}^\infty \mathcal{A}_n \subseteq \mathcal{A}_\perp$ and so  $\underline{\mathcal{A}}\subseteq \mathcal{A}_\perp$.\\
		Now, we show that $\mathcal{A}_\perp \subseteq \overline{\mathcal{A}}$, or equivalently $\mathcal{A}_n \xrightarrow{\perp} \overline{\mathcal{A}}$ (Lemma \ref{lem10}). Suppose that $\mathcal{A}_n \not\xrightarrow{\perp} \overline{\mathcal{A}}$. Then there exists a  sequence $\{A_n : A_n \in \mathcal{A}_n\}_{n \in \mathbb{N}}$ such that $E_{\overline{\mathcal{A}}}^\perp(\chi_{A_n}) \not\xrightarrow{weakly} 0$.  So we can find $g \in L^\psi(\mathcal{A})$ such that
		\begin{align*}
			\int_\Omega g\ E_{\overline{\mathcal{A}}}^\perp(\chi_{A_n}) d \mu \not\rightarrow 0.  \tag{1}
		\end{align*}
		Clearly $\{\chi_{A_{n}}\}_{n \in \mathbb{N}}$ is a bounded sequence in $L^\varphi(\mathcal{A})$, hence there exists a subsequence of $\{\chi_{A_{n_k}}\}_{k \in \mathbb{N}}$ that converges weakly to some $h \in L^\varphi(\mathcal{A})$. Since for all $m \in \mathbb{N}$ and $n_k \geq m$, $\chi_{A_{n_k}}$ is $\bigvee_{r=m}^\infty \mathcal{A}_r$-measurable, then $h$ is $\overline{\mathcal{A}}$-measurable and so $E\left(h \mid \overline{\mathcal{A}}\right) = h$ and 
		\begin{align*}
			\lim_k \int_\Omega g\ (\chi_{A_{n_k}} - E(\chi_{A_{n_k}} \mid \overline{\mathcal{A}}))  d\mu &= \lim_k \int_\Omega E^\perp_{\overline{\mathcal{A}}}(g)  \chi_{A_{n_k}} d\mu\\
			&= \int_\Omega E^\perp_{\overline{\mathcal{A}}}(g) h d\mu\\
			&= \int_\Omega g\ E^\perp_{\overline{\mathcal{A}}}(h) d\mu\\
			& = \int_\Omega g (h - h)d\mu\\
			& = 0,
		\end{align*}
		which contradicts (1). Thus $\mathcal{A}_n \xrightarrow{\perp} \overline{\mathcal{A}}$. 
	\end{proof}
\end{lem}
In the next Lemma, we show that the weakly convergence of the sequence $E_{\mathfrak{D}}^\perp(E\left(f \mid \mathcal{A}_n\right))$, for  all $f \in L^\varphi(\mathcal{A})$, to zero in Orlicz space $L^\varphi(\mathcal{A})$ is a necessary condition for Orthogonally convergence of the sequence $\{\mathcal{A}_n\}_{n \in \mathbb{N}}$ to $\mathfrak{D}$.

\begin{lem}\label{lem12}
	Let $(\Omega, \mathcal{A}, \mu)$ be a measure space, $\varphi\in\Delta_2$ be a Young function,   $\{\mathcal{A}_n\}_{n \in \mathbb{N}}$ be a sequence of $\sigma$-subalgebras of $\mathcal{A}$ and also $\mathfrak{D}$  be a $\sigma$-subalgebra of $\mathcal{A}$ such that   $\mathcal{A}_{n}\xrightarrow{\perp}\mathfrak{D}$. Then $E_{\mathfrak{D}}^\perp(E\left(f \mid \mathcal{A}_n\right)) \xrightarrow{weakly} 0$, for all $f \in L^\varphi(\mathcal{A})$. 
	\begin{proof}
		Let $f \in L^\varphi(\mathcal{A})$ with $0 \leq f \leq 1$. For every $n, N \in \mathbb{N}$, set
		$$A_{n,k} := \left\{\omega \in \Omega : \frac{k-1}{N} \leq E\left(f \mid \mathcal{A}_n\right)(\omega) < \frac{k}{N} \right\},$$
		where $1 \leq k \leq N+1$. It is clear that $A_{n,k} = E\left(f\mid \mathcal{A}_n\right)^{-1}\left[\frac{k-1}{N}, \frac{k}{N}\right)$, and for each  $n \in \mathbb{N}$, the sets $A_{n,k}$ are disjoint and their union is $\Omega$. Now we define 
		$$g_{n,N} := \sum_{k=1}^{N+1} \frac{k-1}{N}\ \chi_{A_{n,k}}.$$
		Hence $E\left(f\mid \mathcal{A}_n\right) = \sum_{k=1}^{N+1} E\left(f \mid \mathcal{A}_n\right) \chi_{A_{n,k}}$, and so 
		\begin{align*}
			|g_{n,N} - E\left(f\mid \mathcal{A}_n\right)| &= \left| \sum_{k=1}^{N+1} \left(\frac{k-1}{N} - E\left(f\mid \mathcal{A}_n\right)\right) \chi_{A_{n,k}} \right|\\
			&\leq \sum_{k=1}^{N+1} \left|\left(\frac{k-1}{N} - \frac{k}{N}\right) \chi_{A_{n,k}}\right| \\
			&= \sum_{k=1}^{N+1} \frac{1}{N} \  \chi_{A_{n,k}}.
		\end{align*}
		Thus 
		\begin{align*}
			\int_\Omega \varphi\left(\frac{N\  |g_{n,N} - E\left(f \mid \mathcal{A}_n\right)|}{\varphi(1)  \mu(\Omega) + 1}\right) d\mu&\leq \frac{1}{\varphi(1) \mu(\Omega) + 1} \int_\Omega \varphi(N\ |g_{n,N} - E\left(f \mid \mathcal{A}_n\right)|) d\mu\\
			&\leq \frac{1}{\varphi(1) \mu(\Omega) + 1} \int_\Omega \varphi\left(N \sum_{k=1}^{N+1} \frac{1}{N} \chi_{A_{n,k}}\right) d\mu \\
			&= \frac{1}{\varphi(1) \mu(\Omega) + 1} \sum_{k=1}^{N+1} \int_\Omega \varphi(\chi_{A_{n,k}}) d\mu\\
			&= \frac{1}{\varphi(1) \mu(\Omega) + 1} \sum_{k=1}^{N+1} \varphi(1)\  \mu(A_{n,k})\\
			& = \frac{\varphi(1)}{\varphi(1) \mu(\Omega) + 1} \sum_{k=1}^{N+1} \mu(A_{n,k})\\
			&= \frac{\varphi(1) \mu(\Omega)}{\varphi(1) \mu(\Omega) + 1} \\
			&< 1,
		\end{align*}
		and consequently
		$$N_\varphi(|g_{n,N} - E\left(f\mid \mathcal{A}_n\right)|) \leq \frac{\varphi(1) \mu(\Omega) + 1}{N}.$$
		The above inequality implies that $\lim_{N \rightarrow \infty} N_\varphi(|g_{n,N} - E\left(f\mid \mathcal{A}_n\right)|) = 0$.  Now, We prove that  $E_{\mathfrak{D}}^\perp(g_{n, N}) \xrightarrow{weakly} 0$. To this end we write 
		\begin{align*}
			E_{\mathfrak{D}}^{\perp}(g_{n, N})&=g_{n, N}-E_{\mathfrak{D}}(g_{n, N})\\
			&=\sum_{k=1}^{N+1}\frac{k-1}{N}\chi_{A_{n,k}}-E_{\mathfrak{D}}(\sum_{k=1}^{N+1}\frac{k-1}{N}\chi_{A_{n,k}})\\
			&=\sum_{k=1}^{N+1}\frac{k-1}{N}(\chi_{A_{n,k}}-E_{\mathfrak{D}}(\chi_{A_{n,k}}))\\
			&=\sum_{k=1}^{N+1}\frac{k-1}{N}E_{\mathfrak{D}}^{\perp}(\chi_{A_{n,k}}).
		\end{align*} 
		As  $\mathcal{A}_{n}\xrightarrow{\perp}\mathfrak{D}$, we have  $E_{\mathfrak{D}}^{\perp}(\chi_{A_{n,k}})\xrightarrow{weakly}0$, and so   $E_{\mathfrak{D}}^{\perp}(g_{n, N})\xrightarrow{weakly}0$, as $n \rightarrow \infty$. 
		Therefore, for every $h\in L^{\psi}(\mathcal{A})$, we have  
		\begin{flalign*}
			\lim_{n}\int_{\Omega}|h  E^{\perp}_{\mathfrak{D}}&(E\left({f}\mid{\mathcal{A}_n}\right))| d\mu\\ &=\lim_{n}\int_{\Omega}|h  E^{\perp}_{\mathfrak{D}}(E\left({f}\mid{\mathcal{A}_n}\right))-h  E^{\perp}_{\mathfrak{D}}(g_{n,N})+h E^{\perp}_{\mathfrak{D}}(g_{n,N})| d\mu\\
			&\le \lim_{n}\int_{\Omega}(|h E^{\perp}_{\mathfrak{D}}(E\left({f}\mid{\mathcal{A}_n}\right)-g_{n,N})|+h E^{\perp}_{\mathfrak{D}}(g_{n,N})|)d\mu\\
			&\le 2\ \lim_{n} \|h\|_{\psi} \
			\|E^{\perp}_{\mathfrak{D}}(E\left({f}\mid{\mathcal{A}_n}\right)-g_{n,N})\|_{\varphi}\\
			&\leq 4 \ \|h\|_{\psi} \  \lim_{n} \|E\left({f}\mid{\mathcal{A}_n}\right)-g_{n,N}\|_{\varphi}\\
			&\leq 4\ \|h\|_{\psi}\  \frac{\varphi(1)\mu(\Omega)+1}{N}.
		\end{flalign*} 
		Since $N$ is arbitrary, then  we have  $E_{\mathfrak{D}}^\perp(E\left(f \mid \mathcal{A}_n\right)) \xrightarrow{weakly} 0$, as $n \rightarrow \infty$, for all $f \in L^\varphi(\mathcal{A})$,  with $0 \leq f \leq 1$. Let $f \in L^\infty(\mathcal{A})$ and  $\tilde{f} := \frac{f + \|f\|_\infty}{2 \|f\|_\infty}$. Clearly $0 \leq \tilde{f} \leq 1$, and  $f = 2 \|f\|_\infty \tilde{f} - \|f\|_\infty$. Using the fact that $E_{\mathfrak{D}}^{\perp}(1)=0$,  linearity of the conditional expectation and its orthogonal complement, we have 
		\begin{align*}
			E_{\mathfrak{D}}^\perp(E\left(f\mid \mathcal{A}_n\right)) &= E_{\mathfrak{D}}^\perp\left(E\left(2 \|f\|_\infty \tilde{f} - \|f\|_\infty \mid \mathcal{A}_n\right)\right)\\
			&= 2 \|f\|_\infty \  E_{\mathfrak{D}}^\perp\left(E\left(\tilde{f} \mid \mathcal{A}_n\right)\right) - \|f\|_\infty E_{\mathfrak{D}}^\perp\left(E\left(1\mid \mathcal{A}_n\right)\right)\\
			&= 2 \|f\|_\infty E_{\mathfrak{D}}^\perp\left(E\left(\tilde{f}\mid \mathcal{A}_n\right)\right) - \|f\|_\infty E_{\mathfrak{D}}^\perp(1)\\
			&= 2 \|f\|_\infty E_{\mathfrak{D}}^\perp\left(E(\tilde{f}\mid \mathcal{A}_n)\right)\\
			& \xrightarrow{weakly} 0.
		\end{align*}
		Since $(\Omega, \mathcal{A}, \mu)$ is a finite measure space, then   $L^\infty(\mathcal{A})$ contains simple functions  and so  $L^\infty(\mathcal{A})$ is a dense subspace of $L^\varphi(\mathcal{A})$. Let $g \in L^\psi(\mathcal{A})$, $f \in L^\varphi(\mathcal{A})$  and $\varepsilon > 0$. Then  there exists $f^{\prime} \in L^\infty(\mathcal{A})$ such that
		$$\|f - f^{\prime}\|_\varphi < \frac{\varepsilon}{2 \|g\|_\psi},$$
		and so 
		\begin{flalign*}
			\big|\int_\Omega g  E_{\mathfrak{D}}^\perp&\left(E\left(f\mid \mathcal{A}_n\right)\right)  d\mu\big| \\
			&\leq \left|\int_\Omega g E_{\mathfrak{D}}^\perp\left(E\left(f - f^{\prime}\mid \mathcal{A}_n\right)\right) d\mu\right| + \left|\int_\Omega g E_{\mathfrak{D}}^\perp(E\left(f^{\prime} \mid \mathcal{A}_n\right)) d\mu\right|\\
			&\leq \|g\|_\psi \|E\left(f - f^{\prime} \mid \mathcal{A}_n\right)\|_\varphi + \left|\int_\Omega g E_{\mathfrak{D}}^\perp(E\left(f^{\prime} \mid \mathcal{A}_n\right))  d\mu\right|\\
			& \leq \|g\|_\psi \|f - f^{\prime}\|_\varphi + \left|\int_\Omega g E_{\mathfrak{D}}^\perp(E\left(f^{\prime}\mid \mathcal{A}_n\right)) d\mu\right|\\
			& < \|g\|_\psi \frac{\varepsilon}{2 \|g\|_\psi} + \left|\int_\Omega g E_{\mathfrak{D}}^\perp(E\left(f^{\prime}\mid \mathcal{A}_n\right)) d\mu\right|\\
			&= \frac{\varepsilon}{2} + \left|\int_\Omega g E_{\mathfrak{D}}^\perp(E\left(f^{\prime}\mid \mathcal{A}_n\right)) d\mu\right|.
		\end{flalign*}
		Since $f^{\prime} \in L^\infty(\mathcal{A})$, by the previous step, $\lim_{n \rightarrow \infty} \int_\Omega g E_{\mathfrak{D}}^\perp(E\left(f^\prime\mid \mathcal{A}_n\right)) d\mu = 0$. Thus,
		$$\lim_{n \rightarrow \infty} \left|\int_\Omega g E_{\mathfrak{D}}^\perp(E\left(f\mid \mathcal{A}_n\right)) d\mu\right| \leq \frac{\varepsilon}{2} + 0 = \frac{\varepsilon}{2}.$$
		Since $\varepsilon > 0$ is arbitrary, we have $\lim_{n \rightarrow \infty} \int_\Omega g\ E_{\mathfrak{D}}^\perp(E\left(f\mid\mathcal{A}_n\right)) d\mu= 0$, for all $g \in L^\psi(\mathcal{A})$. This completes the proof. 
	\end{proof}
\end{lem}

Here a provide a technical lemma for latter use.
\begin{lem}\label{lem13}
	If $(\Omega, \mathcal{A}, \mu)$ is a measure space, $\varphi$ be a Young function and $f \in L^{\varphi}(\mathcal{A})$, then for every  $\sigma$-subalgebras $\mathcal{B}$ and $\mathcal{C}$ of $\mathcal{A}$, we have
	$$\int_D E_{\mathcal{B}}(E\left(f \mid \mathcal{C}\right)) d\mu = \int_D E\left(f \mid \mathcal{B} \cap \mathcal{C}\right) d\mu,$$
	for all $D \in \mathcal{B} \cap \mathcal{C}$.
	\begin{proof}
		Let $D \in \mathcal{B} \cap \mathcal{C}$. Then by definition of conditional expectation, we have 
		\begin{align*}
			\int_D E\left(f\mid \mathcal{B} \cap \mathcal{C}\right) d\mu= \int_D f d\mu. \quad \tag{1}
		\end{align*}
		On the other hand, since $D \in \mathcal{C}$ and $D \in \mathcal{B}$, then we have 
		$$\int_D E\left(f\mid \mathcal{C}\right) d\mu = \int_D f d\mu, $$
		and 
		$$\int_D E_{\mathcal{B}}(E\left(f\mid \mathcal{C}\right)) d\mu = \int_D E\left(f\mid \mathcal{C}\right) d\mu.$$
		By the above observations we get that 
		$$\int_D E_{\mathcal{B}}(E\left(f\mid \mathcal{C}\right)) d\mu = \int_D f d\mu,$$
		and so by (1) we have 
		$$\int_D E_{\mathcal{B}}(E\left(f\mid \mathcal{C}\right)) d\mu = \int_D E\left(f\mid \mathcal{B} \cap \mathcal{C}\right) d\mu,$$
		for every $D \in \mathcal{B} \cap \mathcal{C}$. This complete the prove.  
	\end{proof}
\end{lem}

\begin{lem}
	Let $(\Omega, \mathcal{A}, \mu)$ be a measure space, $\varphi\in\Delta_2$ be a Young function  and $\{\mathcal{A}_n\}_{n \in \mathbb{N}}$ be a sequence of $\sigma$-subalgebras of $\mathcal{A}$. If $\mathfrak{D}$ is also a $\sigma$-subalgebra of $\mathcal{A}$ such that $\mathcal{A}_n \xrightarrow{\perp} \mathfrak{D}$, then we have  $E\left(f\mid \mathcal{A}_n\right) \xrightarrow{L^\varphi(\mathcal{A})} 0$, for all $f \in L^\varphi(\Omega, \mathfrak{D}, \mu)^\perp$.
	\begin{proof}
		Let $f \in L^\varphi(\mathfrak{D})^\perp$. Then, for every $n \in \mathbb{N}$, we get that 
		\begin{align*}
			\|E\left(f\mid \mathcal{A}_n\right)\|_\varphi =& \|E_{\mathfrak{D}}^\perp(E\left(f\mid \mathcal{A}_n\right)) + E_{\mathfrak{D}}(E\left(f \mid \mathcal{A}_n\right))\|_\varphi\\
			\leq& \|E_{\mathfrak{D}}^\perp(E\left(f\mid \mathcal{A}_n\right))\|_\varphi + \|E_{\mathfrak{D}}(E\left(f\mid \mathcal{A}_n\right))\|_\varphi. \tag{1}
		\end{align*}
		By Lemma \ref{lem12} and Jensen inequality, we have 
		\begin{align*}
			\|E_{\mathfrak{D}}^\perp(E\left(f\mid \mathcal{A}_n\right))\|_\varphi &= \inf \left\{ \lambda > 0 : \int_\Omega \varphi\left(\frac{E_{\mathfrak{D}}^\perp(E\left(f \mid \mathcal{A}_n\right))}{\lambda}\right) d\mu \leq 1 \right\}\\ 
			&= \inf \left\{ \lambda > 0 : \frac{1}{\mu(\Omega)} \int_\Omega \varphi\left(\frac{E_{\mathfrak{D}}^\perp(E\left(f \mid \mathcal{A}_n\right))}{\lambda}\right) d\mu \leq \frac{1}{\mu(\Omega)} \right\}\\
			&\leq \inf \left\{ \lambda > 0 : \varphi\left(\frac{1}{\mu(\Omega)} \int_\Omega \frac{E_{\mathfrak{D}}^\perp(E\left(f\mid \mathcal{A}_n\right))}{\lambda} d\mu \right) \leq \frac{1}{\mu(\Omega)} \right\}\\
			&= \inf \left\{ \lambda > 0 : \lambda \geq \frac{\int_\Omega E_{\mathfrak{D}}^\perp(E\left(f\mid \mathcal{A}_n\right))   d\mu}{\mu(\Omega)\ \varphi^{-1} (\frac{1}{\mu(\Omega)})} \right\}\\
			&= \frac{\int_\Omega E_{\mathfrak{D}}^\perp(E\left(f\mid \mathcal{A}_n\right)) d\mu}{\mu(\Omega)\  \varphi^{-1}\left(\frac{1}{\mu(\Omega)}\right)} \\
			&= \frac{\int_\Omega E\left(f\mid \mathcal{A}_n \cap \mathfrak{D}\right) 1_{\mathfrak{D}} d\mu}{\mu(\Omega)\ \varphi^{-1}\left(\frac{1}{\mu(\Omega)}\right)}\\
			& \xrightarrow{n \rightarrow \infty} 0. \tag{2}
		\end{align*}
		Similarly, we have 
		$$\|E_{\mathfrak{D}}(E\left(f \mid\mathcal{A}_n\right))\|_\varphi \leq \frac{\int_\Omega E_{\mathfrak{D}}(E\left(f\mid \mathcal{A}_n\right))  d\mu}{\mu(\Omega)\  \varphi^{-1}\left(\frac{1}{\mu(\Omega)}\right)}.$$
		The Lemma \ref{lem13} implies that 
		$\int_{\Omega}E_{\mathfrak{D}}(E\left({f}\mid{\mathcal{A}_n}\right)) d\mu=\int_{\Omega}E\left({f}\mid{\mathcal{A}_n\cap\mathfrak{D}}\right) d\mu$, for any 
		$n\in\mathbb{N}$. Now, we show that $E\left(f \mid \mathcal{A}_n \cap \mathfrak{D}\right) \stackrel{a.e.}{=} 0$. Since $f \in L^\varphi(\mathfrak{D})^\perp\subseteq L^\psi(\mathfrak{D})$, then for every $g \in L^\varphi(\mathfrak{D})$, we have $\int_\Omega fg d\mu = 0$. Since $\Omega$ is finite, then for every $A \in \mathfrak{D}$,  $\chi_{A}\in L^{\varphi}$ and so $\int_{\Omega} f\chi_{A} d \mu = 0$. Hence $E\left(f\mid \mathfrak{D}\right) \stackrel{a.e.}{=} 0$ and so for every $B \in \mathcal{A}_n \cap \mathfrak{D}$, we have $\int_B E\left(f\mid \mathcal{A}_n \cap \mathfrak{D}\right) d\mu = \int_B f d\mu = 0$. This implies that $E\left(f\mid \mathcal{A}_n \cap \mathfrak{D}\right) \stackrel{a.e.}{=} 0$. Therefore we have 
		\begin{align*}
			\|E_{\mathfrak{D}}(E\left({f}\mid{\mathcal{A}_n}\right))\|_{\varphi}&\leq\frac{\int_{\Omega}E_{\mathfrak{D}}(E\left({f}\mid{\mathcal{A}_n}\right))\ d\mu}{\mu(\Omega)\ \varphi^{-1}(\frac{1}{\mu(\Omega)})}\\
			&=\frac{\int_{\Omega}E\left({f}\mid{\mathcal{A}_n\cap\mathfrak{D}}\right) d\mu}{\mu(\Omega)\ \varphi^{-1}(\frac{1}{\mu(\Omega)})}\\
			&=0.\tag{3}
		\end{align*}
		By applying (2) and (3) in relation (1), we get the proof. 
	\end{proof}
\end{lem}

By the last observations, now in the next Corollary we have an equivalent condition to the orthogonally convergence of a sequence of $\sigma$-subalgebras.
\begin{cor}
	Let $(\Omega, \mathcal{A}, \mu)$ be a measure space, $\varphi\in\Delta_2$ be a Young function and $\{\mathcal{A}_n\}_{n=1}^\infty$ be a sequence of $\sigma$-subalgebras of $\mathcal{A}$. Then  $\mathcal{A}_n \xrightarrow{\perp} \underline{\mathcal{A}}$ if and only if 
	$$E\left(f\mid \mathcal{A}_n\right) \xrightarrow{L^\varphi(\mathcal{A})} E\left(f\mid \underline{\mathcal{A}}\right),$$
	for all $f \in L^\varphi(\mathcal{A})$.
\end{cor}

\section{ \sc\bf $\mu\perp$-Convergence}
In this section we study the sequences of $\sigma$-subalgebras that are convergent in both $\mu$-convergence $\perp$-convergence to a common limit. Let $(\Omega, \mathcal{A}, \mu)$ be a measure space, $\{\mathcal{A}_n\}_{n \in \mathbb{N}}$ be a sequence of $\sigma$-subalgebras of $\mathcal{A}$ and $\mathfrak{D}$ also be a $\sigma$-subalgebra of $\mathcal{A}$. We say that $\{\mathcal{A}_n\}_{n \in \mathbb{N}}$ $\mu\perp$-converges  to $\mathfrak{D}$, and  denote it by $\mathcal{A}_n \xrightarrow{\mu\perp} \mathfrak{D}$, if $\mathcal{A}_n \xrightarrow{\mu} \mathfrak{D}$ and $\mathcal{A}_n \xrightarrow{\perp} \mathfrak{D}$.

\begin{lem}\label{lem4.1}
	Let $(\Omega, \mathcal{A}, \mu)$ be a measure space and $\{\mathcal{A}_n\}_{n \in \mathbb{N}}$ be a sequence of $\sigma$-subalgebras of $\mathcal{A}$. Then 
	$$\underline{\mathcal{A}} \subseteq \mathcal{A}_\mu \subseteq \mathcal{A}_\perp \subseteq \overline{\mathcal{A}}.$$
	\begin{proof}
		As is proved in the Lemmas \ref{lem1}, \ref{lem10} we have  $\underline{A} \subseteq A_\mu \subseteq \overline{A}$ and $\underline{\mathcal{A}}\subseteq\mathcal{A}_\perp \subseteq \overline{\mathcal{A}}$. So we have to show that  $\mathcal{A}_{\mu}\subseteq\mathcal{A}_{\perp}$. To this end, let $A \in \mathcal{A}_{\mu}$. Then, there exists a sequence  $\{A_n: A_n\in\mathcal{A}_n\}_{n \in \mathbb{N}}$  such that $\lim_{n}\mu(A_n \Delta A)=0$. This implies that  $\chi_{A_n} \xrightarrow{\mu} \chi_A$ and so by the Lemma \ref{lem3} we have  $N_{\varphi}(\chi_{A_n} - \chi_A) \rightarrow 0$.  Since strong convergence in any Banach space (including $L^\varphi$), implies weak convergence, we have
		$$\chi_{A_n} \xrightarrow{\text{weakly}} \chi_A.$$ 
		Hence, by the definition of $\mathcal{A}_{\perp}$, the weak limit $\chi_A$ is an element of $W_{\varphi}$. As $\mathcal{A}_{\perp}$ is the smallest $\sigma$-algebra with respect to which all functions in $W_{\varphi}$ are measurable, $\chi_A$ must be $\mathcal{A}_{\perp}$-measurable. Therefore, $A \in \mathcal{A}_{\perp}$. This completes the proof. 
	\end{proof}
\end{lem}

In the next Proposition we characterize $\mu\perp$-Convergence of a sequence of $\sigma$-subalgebras.

\begin{prop}
	Let $(\Omega, \mathcal{A}, \mu)$ be a measure space,  $\{\mathcal{A}_n\}_{n \in \mathbb{N}}$ be a sequence of $\sigma$-subalgebras of $\mathcal{A}$ and  $\mathfrak{D}$ also be a $\sigma$-subalgebra of $\mathcal{A}$.  Then $\mathcal{A}_n \xrightarrow{\mu\perp} \mathfrak{D}$ if and only if $\mathcal{A}_\mu=\mathfrak{D}=\mathcal{A}_\perp$.
	\begin{proof}
		It is a direct consequence of the lemmas \ref{lem1}, \ref{lem10} and \ref{lem4.1}.
	\end{proof}
\end{prop}

In the following Lemma we find a necessary condition for $\mu\perp$--Convergence of a sequence of $\sigma$-subalgebras in the setting of Orlicz spaces.
\begin{lem}\label{lem15}
	Let $(\Omega, \mathcal{A}, \mu)$ be a measure space,  $\varphi\in\Delta_2$ be a Young function, $\{\mathcal{A}_n\}_{n \in \mathbb{N}}$ be a sequence of $\sigma$-subalgebras of $\mathcal{A}$ and  $\mathfrak{D}$ be a $\sigma$-subalgebra of $\mathcal{A}$. If $\mathcal{A}_n \xrightarrow{\mu\perp} \mathfrak{D}$, then 
	$$ E\left(f\mid\mathcal{A}_n\right)-E\left(f\mid\mathfrak{D}\right) \xrightarrow{weakly} 0, $$ 
	for all $f \in L^\varphi(\mathcal{A})$.
	\begin{proof}
		To prove the assertion we need to show that for every $f \in L^\varphi(\mathcal{A})$,   
		$$\lim_{n \rightarrow \infty}\int_\Omega (E\left(f\mid\mathcal{A}_n\right)-E\left(f\mid\mathfrak{D}\right))g d\mu = 0,$$
		for all $g \in L^\psi(\mathcal{A})$. Since $\mathcal{A}_n \xrightarrow{\mu} \mathfrak{D}$, the Lemma \ref{lem3} implies that 
		$$E_{{\mathcal{A}}_n}(E\left({f}\mid{\mathfrak{D}}\right))\xrightarrow{L^{\varphi}(\mathcal{A})}E\left({f}\mid{\mathfrak{D}}\right).$$
		On the other hand, we have $\mathcal{A}_n \xrightarrow{\perp} \mathfrak{D}$ and so the Lemma \ref{lem12} implies that $$E_{\mathfrak{D}}^\perp(E\left(f\mid\mathcal{A}_n\right)) \xrightarrow{weakly} 0.$$  
		Hence
		$$\lim_{n \rightarrow \infty} \int_\Omega E_{\mathfrak{D}}^\perp(E\left(f\mid\mathcal{A}_n\right)) g d\mu=0,$$
		and in particular for $g=1$, we have
		$$\lim_{n \rightarrow \infty} \int_\Omega (E\left(f \mid \mathcal{A}_n\right) - E_{\mathfrak{D}}(E\left(f \mid \mathcal{A}_n\right))) d\mu = 0.$$
		Now by using the Lemma  \ref{lem13}, we have  
		\begin{align*}
			\int_{\Omega}E\left({f}\mid{\mathcal{A}_n}\right) d\mu&-\int_{\Omega}E_{\mathfrak{D}}E\left({f}\mid{\mathcal{A}_n}\right) d\mu\xrightarrow{n}0\\
			&\Rightarrow \int_{\Omega}E\left({f}\mid{\mathcal{A}_n}\right) d\mu-\int_{\Omega}E\left({f}\mid{\mathfrak{D}\cap\mathcal{A}_n}\right) d\mu\xrightarrow{n}0\\
			&\Rightarrow 
			\int_{\Omega}E\left({f}\mid{\mathcal{A}_n}\right) d\mu-\int_{\Omega}E\left({f}\mid{\mathcal{A}_n}\cap\mathfrak{D}\right) d\mu\xrightarrow{n}0\\
			&\Rightarrow 
			\int_{\Omega}E\left({f}\mid{\mathcal{A}_n}\right) d\mu-\int_{\Omega}E_{\mathcal{A}_n}E\left({f}\mid{\mathfrak{D}}\right) d\mu\xrightarrow{n}0, \tag{2}
		\end{align*} 
		and so 
		\begin{flalign*}
			\|&E\left({f}\mid{\mathcal{A}_n}\right)-E_{\mathcal{A}_n}(E\left({f}\mid{\mathfrak{D}}\right))\|_{\varphi}\\
			&=\inf\Big\{\lambda>0: \int_{\Omega}\varphi\Bigl(\frac{E\left({f}\mid{\mathcal{A}_n}\right)-E_{\mathcal{A}_n}(E\left({f}\mid{\mathfrak{D}}\right))}{\lambda}\Bigr) d\mu\leq1\Big\}\\
			&= \inf\Big\{\lambda>0: \frac{1}{\mu(\Omega)}\int_{\Omega}\varphi\Bigl(\frac{E\left({f}\mid{\mathcal{A}_n}\right)-E_{\mathcal{A}_n}(E\left({f}\mid{\mathfrak{D}}\right))}{\lambda}\Bigr)d\mu\leq\frac{1}{\mu(\Omega)}\Big\}\\
			&\leq \inf\Big\{\lambda>0: \varphi\Bigl(\frac{1}{\mu(\Omega)}\int_{\Omega}\frac{E\left({f}\mid{\mathcal{A}_n}\right)-E_{\mathcal{A}_n}(E\left({f}\mid{\mathfrak{D}}\right))}{\lambda}d\mu\Bigr)\leq\frac{1}{\mu(\Omega)}\Big\}\\
			&=\inf\Big\{\lambda>0: \lambda\geq \frac{\int_{\Omega}E\left({f}\mid{\mathcal{A}_n}\right)-E_{\mathcal{A}_n}(E\left({f}\mid{\mathfrak{D}}\right)) d\mu}{\mu(\Omega)\varphi^{-1}(\frac{1}{\mu(\Omega)})}\Bigr\}\\
			&=
			\frac{\int_{\Omega}E\left({f}\mid{\mathcal{A}_n}\right)-E_{\mathcal{A}_n}(E\left({f}\mid{\mathfrak{D}}\right))d\mu}{\mu(\Omega)\varphi^{-1}(\frac{1}{\mu(\Omega)})}
			\\
			&=\frac{\int_{\Omega}E\left({f}\mid{\mathcal{A}_n}\right)d\mu-\int_{\Omega}E_{\mathcal{A}_n}(E\left({f}\mid{\mathfrak{D}}\right))d\mu}{\mu(\Omega)\varphi^{-1}(\frac{1}{\mu(\Omega)})}\\
			&\xrightarrow{(2)}0.
			\tag{3}
		\end{flalign*}
		Now for any $g\in L^{\psi}(\mathcal{A})$, the relations (1) and (3) imply that 
		\begin{align*}
			\int_{\Omega}\Bigl(E\left({f}\mid{\mathcal{A}_n}\right)-E\left({f}\mid{\mathfrak{D}}\right)\Bigr)g d\mu&= \int_{\Omega} \Bigl(E_{\mathcal{A}_n}(E\left({f}\mid{\mathfrak{D}}\right))-E\left({f}\mid{\mathfrak{D}}\right)\Bigr) g  d\mu\\
			& \ \ + \int_{\Omega} \Bigl(E\left({f}\mid{\mathcal{A}_n}\right)-E_{\mathcal{A}_n}(E\left({f}\mid{\mathfrak{D}}\right))\Bigr) g d\mu\\
			&\leq 2\|E_{\mathcal{A}_n}(E\left({f}\mid{\mathfrak{D}}\right))-E\left({f}\mid{\mathfrak{D}}\right)\|_{\varphi}\|g\|_{\psi}\\
			& \ \ +2\|E\left({f}\mid{\mathcal{A}_n}\right)-E_{\mathcal{A}_n}(E\left({f}\mid{\mathfrak{D}}\right))\|_{\varphi}\|g\|_{\psi}\\
			&\xrightarrow{n\rightarrow\infty} 0.
		\end{align*}
		This complete the proof.
	\end{proof}
\end{lem}

Finally in the next theorem we find an equivalent condition to $\mu\perp$-convergence of a sequence of $\sigma$-subalgebras in the setting of Orlicz spaces.

\begin{thm}\label{thm4.5}
	Let $(\Omega, \mathcal{A}, \mu)$ be a measure space, $\varphi\in\Delta_2$  be a Young function,  $\{\mathcal{A}_n\}_{n \in \mathbb{N}}$ be a sequence of $\sigma$-subalgebras of $\mathcal{A}$ and also $\mathfrak{D}$ be a $\sigma$-subalgebra of $\mathcal{A}$. Then $\mathcal{A}_n \xrightarrow{\mu\perp} \mathfrak{D}$ if and only if 
	$$
	E\left(f \mid \mathcal{A}_n\right) \xrightarrow{L^\varphi(\mathcal{A})} E\left(f\mid \mathfrak{D}\right),
	$$
	for all $f \in L^\varphi(\mathcal{A})$.
	\begin{proof}
		$(\Rightarrow)$ Suppose that  $\mathcal{A}_n \xrightarrow{\mu\perp} \mathfrak{D}$ and $f \in L^\varphi(\mathcal{A})$. By Lemma \ref{lem15} we get that $E\left(f \mid \mathcal{A}_n\right) - E\left(f\mid \mathfrak{D}\right) \xrightarrow{weakly} 0$. Hence 
		\begin{align*}
			\|E\left({f}\mid{\mathcal{A}_n}\right)&-E\left({f}\mid{\mathfrak{D}}\right)\|_{\varphi}\\
			&=\inf\Big\{\lambda>0: \int_{\Omega}\varphi\Bigl(\frac{E\left({f}\mid{\mathcal{A}_n}\right)-E\left({f}\mid{\mathfrak{D}}\right)}{\lambda}\Bigr)d\mu\leq1\Big\}\\
			&= \inf\Big\{\lambda>0: \frac{1}{\mu(\Omega)}\int_{\Omega}\varphi\Bigl(\frac{E\left({f}\mid{\mathcal{A}_n}\right)-E\left({f}\mid{\mathfrak{D}}\right)}{\lambda}\Bigr)d\mu\leq\frac{1}{\mu(\Omega)}\Big\}\\
			& \leq \inf\Big\{\lambda>0: \varphi\Bigl(\frac{1}{\mu(\Omega)}\int_{\Omega}\frac{E\left({f}\mid{\mathcal{A}_n}\right)-E\left({f}\mid{\mathfrak{D}}\right)}{\lambda}d\mu\Bigr)\leq\frac{1}{\mu(\Omega)}\Big\}\\
			&=\inf\Big\{\lambda>0: \lambda\geq \frac{\int_{\Omega}\Bigl(E\left({f}\mid{\mathcal{A}_n}\right)-E\left({f}\mid{\mathfrak{D}}\right)\Bigr)d\mu}{\mu(\Omega)\varphi^{-1}(\frac{1}{\mu(\Omega)})}\Bigr\}\\
			&= \frac{\int_{\Omega}\Bigl(E\left({f}\mid{\mathcal{A}_n}\right)-E\left({f}\mid{\mathfrak{D}}\right)\Bigr)d\mu}{\mu(\Omega)\varphi^{-1}(\frac{1}{\mu(\Omega)})}\\
			&=\frac{\int_{\Omega}\Bigl(E\left({f}\mid{\mathcal{A}_n}\right)-E\left({f}\mid{\mathfrak{D}}\right)\Bigr) \ 1_{\mathfrak{D}}  d\mu}{\mu(\Omega)\varphi^{-1}(\frac{1}{\mu(\Omega)})}\\
			& \xrightarrow{n\rightarrow\infty}0.
		\end{align*}
		Thus the forward direction follows.\\ 
		$(\Leftarrow)$ Suppose that for every \( f \in L^\varphi(\mathcal{A})\),
		$$ E\left(f\mid \mathcal{A}_n\right) \xrightarrow{L^\varphi(\mathcal{A})} E\left(f\mid\mathfrak{D}\right). $$
		Let $D\in\mathfrak{D}$ and set 
		$$A_n := \left\{ w \in \Omega : E(f \mid \mathcal{A}_n)(w)> \frac{1}{2} \right\},$$ 
		for all  $n \in \mathbb{N}$.  Similar to the final part $(ii)$ of Lemma \ref{lem3} , we get that  
		$$\lim_n \mu(A_n \triangle D)=0,$$ and so $D\in\mathcal{A}_{\mu}$. Thus  $\mathfrak{D}\subseteq\mathcal{A}_{\mu} $ and by the Lemma \ref{lem1} we have 
		$$\mathcal{A}_n \xrightarrow{\mu} \mathfrak{D}.$$
		Now we show that
		$\mathcal{A}_n \xrightarrow{\perp} \mathfrak{D}$.  To this end, let \( g \in W_\varphi \). Then  there exists the sequence \( \{ A_{n_k} : A_{n_k} \in \mathcal{A}_{n_k} \}_{k=1}^\infty \)  such that
		$$\chi_{A_{n_k}} \xrightarrow{weakly} g.$$
		Hence, for every $ h \in L^\psi(\mathcal{A})$, the followings hold:
		\begin{align*}
			\int_{\Omega}h  E\left({g}\mid{\mathfrak{D}}\right)d\mu &=\int_{\Omega}\  E\left({h}\mid{\mathfrak{D}}\right) g  d\mu \\
			&=\lim_{k} \int_{\Omega}  E\left({h}\mid{\mathfrak{D}}\right) \chi_{A_{n_k}}d\mu \\ 
			&=\lim_{k} \int_{\Omega} h  E\left({\chi_{A_{n_k}}}\mid{\mathfrak{D}}\right) d\mu\\ 
			&=\lim_k \Bigl(\int_{\Omega} h (E\left({\chi_{A_{n_k}}}\mid{\mathfrak{D}}\right)-\chi_{A_{n_k}})d\mu \ +\int_{\Omega} h  \chi_{A_{n_k}}d\mu \Bigr)\\
			&=\lim_{k}\int_{\Omega} h {\chi}_{A_{n_k}}d\mu\\
			&=\int_{\Omega}h g \  d\mu.
		\end{align*} 
		Therefore $ g $ and $ E(g \mid \mathfrak{D}) $ are equal almost everywhere, and so the function $ g $ is $ \mathfrak{D} $-measurable.  
		This means that $ W_\varphi$  consists of $ \mathfrak{D} $-measurable functions. Since $ \mathcal{A}_\perp $ is the smallest $\sigma$-subalgebra such that all members of \( W_\varphi \) are measurable with respect to $ \mathcal{A}_\perp $, we have  
		$$\mathcal{A}_\perp\subseteq \mathfrak{D}, $$
		and so the by the Lemma \ref{lem10} we get that  
		$$\mathcal{A}_n \xrightarrow{\perp} \mathfrak{D}.$$
		This completes the proof. 
	\end{proof}
\end{thm}

\begin{exam}
	Let $(\Omega, \mathcal{A}, \mu) = ([0, 1], \mathcal{B}([0, 1]), m)$ be the measure space with the Lebesgue measure $m$ and $\varphi\in \Delta_2$ be a Young function. 
	For each $n \in \mathbb{N}$, let $\mathcal{G}_n$ be the $\sigma$-subalgebra generated by the partition $P_n = \{I_k^{(n)}\}_{k=1}^n$ of $[0, 1]$, where $I_k^{(n)} = \left[\frac{k-1}{n}, \frac{k}{n}\right)$.
	The Lebesgue measure of each partition element is $m(I_k^{(n)}) = \frac{1}{n}$.
	The conditional expectation $E(f \mid \mathcal{G}_n)$, for any $f \in L^{\varphi}([0, 1])$, is the $\mathcal{G}_n$-measurable function defined as 
	\begin{align*}
		E(f \mid \mathcal{G}_n)(x) &= \sum_{k=1}^{n} \left( \frac{1}{m(I_k^{(n)})} \int_{I_k^{(n)}} f(t) dt \right) \mathbf{1}_{I_k^{(n)}}(x)\\
		& = \sum_{k=1}^{n} \left( n \int_{I_k^{(n)}} f(t) dt \right) \mathbf{1}_{I_k^{(n)}}(x).
	\end{align*}
	Set  
	$$
	\mathfrak{D}:= \sigma\left(\bigcup_{n=1}^\infty \mathcal{G}_n\right).
	$$
	The union $\bigcup_{n=1}^\infty \mathcal{G}_n$ contains all intervals of the form $\left[\frac{k}{n}, \frac{k+1}{n}\right)$, for any $n \in \mathbb{N}$ and $1 \leq k \leq n-1$. Since the algebra of finite unions of such dyadic intervals is dense in the Borel $\sigma$-algebra $\mathcal{A}$, we have 
	$$
	\mathfrak{D} = \mathcal{B}([0, 1]) = \mathcal{A}.
	$$
	Thus, the condition for convergence of the $\sigma$-algebras is satisfied  $\mathcal{G}_n \xrightarrow{m\perp} \mathcal{A}$ and the theorem \ref{thm4.5} implies that 
	$$
	\lim_{n \to \infty} E\left(f \mid \mathcal{G}_n\right) = E\left(f\mid \mathfrak{D}\right).
	$$
\end{exam}

\end{document}